\documentclass[preprint,12pt]{elsarticle}

\usepackage{amsmath,amsfonts,amsthm,amsopn,color,amssymb,graphics,graphicx}
\usepackage{graphicx}
\usepackage{epsfig}
\usepackage{subfigure}
\usepackage{pstricks}
\usepackage{pstricks-add}
\usepackage{amsmath}
\usepackage{amssymb}
\usepackage{amsfonts}
\usepackage{color}
\usepackage{amscd}
\usepackage{latexsym}
\usepackage{mathrsfs}
\usepackage{pst-all}
\usepackage{graphicx,graphics}
\usepackage{subfigure}
\usepackage{srcltx}
\usepackage{psfrag}
\usepackage{graphicx}
\usepackage{epstopdf}
\usepackage{multirow}

\renewcommand{\O}{\Omega}

\newtheorem{theorem}{Theorem}[section]
\newtheorem{lemma}[theorem]{Lemma}

\newtheorem{remark}{Remark}[section]
\newtheorem{example}{Example}

\numberwithin{equation}{section}
\def\<{\leftangle}
\def\>{\rightangle}
\def\({\left(}
\def\){\right)}

\def\<{\langle}
\def\>{\rangle}
\def\q{\quad}

\def\a{\alpha}

\def\g{\gamma}

\def\d{\delta}
\def\h{\hbox}

\def\l{\lambda}

\def\w.r.t.{with respect to}

\def\O{\Omega}

\def\ds{\displaystyle}

\begin{document}

\begin{frontmatter}
\title{Discretized Tikhonov regularization for Robin boundaries localization}
\author[Cao]{Hui Cao}
\address[Cao]{School of Mathematics and Computational Science, Sun Yat-sen University, Guangzhou, 510275, P.R. China}
\ead[Cao]{caohui6@mail.sysu.edu.cn}

\author[Pereverzev]{Sergei V. Pereverzev}
\address[Pereverzev]{Johann Radon Institute for Computational and Applied Mathematics, Austrian Academy and Sciences, Linz, A-4040, Austria}
\ead[Pereverzev]{sergei.pereverzyev@oeaw.ac.at}

\author[Sincich]{Eva Sincich}
\address[Sincich]{Laboratory for Multiphase Processes, University of Nova Gorica, Slovenia}
\ead[Sincich]{eva.sincich@ung.si}

\date{}
\begin{abstract}
We deal with a boundary detection problem arising in nondestructive testing of materials. The problem consists in recovering an unknown portion of the boundary, where a Robin condition is satisfied, with the use of a Cauchy data pair collected on the accessible part of the boundary. We combine a linearization argument with a Tikhonov regularization approach for the local reconstruction of the unknown defect. Moreover, we discuss the regularization parameter choice by means of the so called balancing principle and we present some numerical tests that show the efficiency of our method.
\end{abstract}

\end{frontmatter}

\section{Introduction}

In this paper we deal with an inverse problem arising in corrosion detection. We consider a domain $\Omega\subset \mathbb{R}^2$ which models a $2D$ transverse section of a thin metallic specimen whose boundary is partly accessible and stays in contact with an aggressive environment. Hence, in order to detect the damage which is expected to occur in such a portion of the boundary, one has to rely on the electrostatic measurements of a potential $u$ performed on the accessible portion.

We are then lead to the study of the following elliptic boundary value problem
\begin{equation}\label{P}
\left\{
\begin{array}
{lcl}
\Delta u=0\ ,& \mbox{in $\Omega$ ,}
\\
\dfrac{\partial u}{\partial \nu} =\Phi \ , & \mbox{on $\Gamma_A$ ,}
\\
\dfrac{\partial u}{\partial \nu} + \gamma u=0 \ , & \mbox{on $\Gamma_I$ ,}
\\
u=0 \ , & \mbox{on $\Gamma_D$ .}
\end{array}
\right.
\end{equation}

According to this model $u$ is the harmonic potential in $\Omega$. We assume that the boundary of $\Omega$ is decomposed in three open and disjoint subsets $\Gamma_A, \Gamma_I, \Gamma_D$. On the portion $\Gamma_A$, which is the one accessible to direct inspection, we prescribe a current density $\Phi$  and we measure the corresponding voltage potential $u|_{\Gamma_A}$. The portion $\Gamma_I$, where the corrosion took place, is out of reach. On such a portion the potential $u$ satisfies an homogeneous Robin condition, which models a resistive coupling with the exterior environment by means of the impedance coefficient $\gamma$.

In this paper we are interested in the numerical reconstruction issue of the unknown and damaged boundary $\Gamma_I$ from the data collected on the accessible part of the boundary $\Gamma_A$, that is the Cauchy data pair $(u|_{\Gamma_A}, \Phi)$.

Boundary and parameter identification results related to this stationary inverse problem has been provided by many authors \cite{ADR, AS1, AS2, CFS, CaPS, B, CakK, ChJ, FI, FIM, KaSa, KuS, S1, S2}.

Local uniqueness and conditional stability results for the inverse problem at hand are contained in \cite{CFS} and constitute the theoretical setting on which our numerical analysis relies.  The present local determination of corroded boundaries consists in the localization of a small perturbation $\Gamma_{I, \theta}$ of a reference boundary $\Gamma_I$. It is convenient to introduce a small vector field $\theta\in C^{1}_0(\Gamma_I)$ so that the damaged domain $\Omega_{\theta}$ is such that
\begin{eqnarray*}
\partial \Omega_{\theta}= \overline{\Gamma_A}\cup\overline{\Gamma_D}\cup\overline{\Gamma_{I,\theta}}, 	 	 \ \ \ \Gamma_{I,\theta}=\{z\in \mathbb{R}^2: z= w+ \theta(w),\ w\in \Gamma_I \}.
\end{eqnarray*}
Such a local approach combined with a linearization argument (see \cite{CFS}) allows a reformulation of the problem of the localization of the unknown defect $\Gamma_{I,\theta}$ as the identification of the unknown term $\theta$ in a boundary condition of the type
\begin{eqnarray*}
\dfrac{\partial u'}{\partial \nu} + \gamma u'= \frac{\mbox{d}}{\mbox{d}s}\left(\theta\cdot\nu \frac{\mbox{d}} {\mbox{d}s}u\right) + \gamma \theta\cdot\nu (\gamma + 2H)u
\end{eqnarray*}
at the portion $\Gamma_{I}$, where $u'$ is a harmonic function satisfying homogeneous Neumann and Dirichlet conditions on $\Gamma_A$ and $\Gamma_D$ respectively, $u$ is the solution of (\ref{P}), and $H$ denotes the  curvature of the reference boundary $\Gamma_I$.
As in \cite{CFS} we carry over our analysis under the a-priori assumption of a constant $\gamma$ such that $2H(x)+ \gamma >0$ in $\Gamma_I$ and we limit ourselves to the case of positive fluxes $\Phi$ only.

We linearize the forward map $F: \theta\mapsto u_{\theta}|_{\Gamma_{A}}$, where $u_{\theta}$ is the solution of the system (\ref{Ptheta}) below,  by its Fr\'{e}chet derivative $F'$ and take the \emph{voltage contrast} on $\Gamma_{A}$, as the noisy right-hand term for the considered operator equation,
 $$
 F'\theta=(u_{\theta}-u)|_{\Gamma_{A}}.
 $$
 As in \cite{FIM}, we assume that the unavoidable measurement errors in  \emph{voltage contrast} are not smaller than the truncation error, $o(\|\theta\|)$. Therefore,  if the noise level for  voltage measurements is assumed to be $\d$, then the noise level for the right-hand term of the above operator equation can be written as $\tilde{\d}=K\d$, where a constant $K$ is not necessary to be precisely known.
Our method is based on a discretized Tikhonov regularization argument  where the regularization parameter is chosen by a balancing principle (cf.\cite{CaP1,CaP2,Lu}). Such an \emph{a posteriori} parameter choice can lead to a regularized solution with order-optimal accuracy. At the same time it can provide a reliable estimate for the constant $K$.

\section{Local identification of the unknown boundary}

In this section we shall collect the main identifiability results of which our reconstruction procedure and our numerical tests are a follow up. For a more detailed description we refer to \cite{CFS}.

We denote with $\nu$ the outward normal to $\Gamma_I$ and we assume that $\theta$ is a vector field in $C^1_0(\Gamma_I)$ having a nontrivial normal component $\theta_{\nu}$ on $\Gamma_I$.

Let the Sobolev space $H^1_0(\Omega, \Gamma_D)$ be defined as follows
\begin{eqnarray}
H^1_0(\Omega, \Gamma_D)=\{v\in H^1(\Omega)\ : \ v=0 \ \mbox{on}\ \Gamma_D \ \mbox{in the trace sense} \} \ .
\end{eqnarray}

We introduce the forward map $F$
\begin{eqnarray}\label{fm}
   F : \quad C^1_0(\Gamma_I)
   & \rightarrow & H^{\frac{1}{2}}(\Gamma_A) \\
   \theta & \mapsto &
   u_{\theta}|_{\Gamma_A} \nonumber
\end{eqnarray}
where $u_{\theta}\in H^1_0(\Omega, \Gamma_D)$ is the solution to the elliptic problem
\begin{equation}\label{Ptheta}
   \left\{ \begin{array}{ll}
   \Delta u_{\theta}=0\ & \mbox{in $\Omega_{\theta}$ } \\
   \dfrac{\partial u_{\theta}}{\partial \nu}=\Phi\ &
   \mbox{on $\Gamma_A$ }  \\
   \dfrac{\partial u_{\theta}}{\partial \nu} +
   \gamma u_{\theta} =0  & \mbox{on $\Gamma_{I, \theta}$ }  \\
   u=0 & \mbox{on $\Gamma_{D}$.}
   \end{array} \right.
\end{equation}

We recall the following differentiability property for the forward map $F$.
\begin{lemma}\label{Fderivative}
The operator $F$ in \eqref{fm} is Fr\'{e}chet differentiable
at zero. Indeed,
consider the linear operator
$F': C^1_0(\Gamma_I) \rightarrow H^{\frac{1}{2}}(\Gamma_A)$
defined as
$F'\theta=u'|_{\Gamma_A}$, where $u'$ is the solution to the
boundary value problem
\begin{equation}\label{uprime}
   \left\{\begin{array}{ll}
   \Delta u'=0 & \mbox{in $\Omega$} \\
   \dfrac{\partial u'}{\partial \nu}=0   & \mbox{on $\Gamma_A$} \\
   \dfrac{\partial u'}{\partial \nu} + \gamma u' =
   \displaystyle{\frac{\mathrm{d}}{\mathrm{d}s}
   \left({\theta}_{\nu}\frac{\mathrm{d}}{\mathrm{d}s}u \right) +
   \gamma{\theta}_{\nu}
   \left(\gamma  + 2H \right)u}  & \mbox{on $\Gamma_I$} \\
   u'=0 & \mbox{on $\Gamma_{D}$,}
\end{array}
\right.
\end{equation}
the function $u$ is the solution of \eqref{P}
and $H$ denotes the  curvature of the boundary $\Gamma_I$.
Then,
$$
   \frac{1}{\|\theta\|_{C_0^1(\Gamma_I)}}
   \| F(\theta)-F(0) -
   F'\theta\|_{H^{\frac{1}{2}}(\Gamma_A)}\rightarrow 0
   \qquad \mbox{as $\theta \rightarrow 0$ in $C^1_0(\Gamma_I)$} .
$$
\end{lemma}

Let us also recall that a weak solution to \eqref{uprime} is a function $u'\in H^1_0(\Omega, \Gamma_D)$ such that
\begin{eqnarray}\label{wf}
\int_{\Omega}\nabla u'\nabla v + \int_{\Gamma_I}\gamma u'v = \int_{\Gamma_I}\gamma \theta_{\nu}(\gamma + H)uv - \int_{\Gamma_I}\theta_{\nu}\frac{\mathrm{d}}{\mathrm{d}s}u\frac{\mathrm{d}}{\mathrm{d}s}v
\end{eqnarray}
for all $v \in H^1_0(\Omega, \Gamma_D). $

The following theorem ensures that the operator $F'$ is injective, under some reasonable hypothesis. This property allows us to conclude that the solution $\theta$ to our inverse problem is identifiable, at least for small perturbations. Moreover,  we recall a conditional Lipschitz type upper bound for $\theta$
on a suitable portion of $\Gamma_I$
in terms of $u'|_{\Gamma_A} = F'\theta$,
thus showing that the inversion of $F'$
is not too much ill-behaved.
\begin{theorem}\label{iniettivo}
Let $\Phi\in H^{\frac{1}{2}}(\Gamma_A)$ be nonnegative in the sense of distributions.
Let us assume that $2H(x) + \gamma > 0 $ and
$\theta_{\nu}(x)\leqslant 0$ for any $x\in \Gamma_I$.
Then $F'$ is injective. Moreover, there exists a positive constant ${c}>0$  such that
$$
   \|u'\|_{H^{\frac{1}{2}}(\Gamma_A)}\geqslant
{c}\int_{\tilde{\Gamma}_I}|\theta|,
$$
where $\tilde{\Gamma}_I$ is an inner portion of the boundary $\Gamma_I$.
\end{theorem}

Finally, in the next theorem, we consider $L^{2}(\Gamma_A)$
as codomain space of the operator $F'$ introduced
in Lemma \ref{Fderivative} stating a compactness result.
\begin{theorem}   \label{teo:compact}
The linear operator
\begin{eqnarray*}
   F': \quad C^1_0(\Gamma_I)
   &\rightarrow & L^{2}(\Gamma_A) \nonumber \\
   \theta & \mapsto & u'|_{\Gamma_A}
\end{eqnarray*}
where $u'$ is the solution to the
boundary value problem \eqref{uprime}, is compact.
\end{theorem}

\section{Tikhonov regularization for a local reconstruction and an estimate of the accuracy}

Here and in the following, with a slight abuse of notation, we shall denote by $F'$ the compact operator
\begin{eqnarray}\label{hilbert}
F': && H^1_0(\Gamma_I) \rightarrow L^2(\Gamma_A)\\
&&\ \ \ \ \  \theta \ \ \ \mapsto \ \  u', \nonumber
\end{eqnarray}
where $u'$ satisfies the weak formulation in \eqref{wf} for any $v\in H^1_0(\Omega, \Gamma_D)$.

The existence and uniqueness of $u'\in H^1_0(\Omega,\Gamma_D)$ follows from standard arguments on elliptic boundary value problems. Moreover, the compactness of $F'$ in \eqref{hilbert} follows along the lines of Theorem 4.5 in \cite{CFS}.

In view of this compactness property, the issue of the identification of $\theta$
may be interpreted as the regularized inversion of the above compact operator $F'=F'(\Gamma_{I})$ between the Hilbert spaces $H^1_0(\Gamma_I)$ and $L^2(\Gamma_A)$. Such kind of
reformulation allows us to deal with the approximate inversion by the technique of Tikhonov regularization.

 We are interested in finding the  solution to operator equation
\begin{equation}\label{3.1}
  F' \theta =\bar{r}:=u'|_{\Gamma_A},
\end{equation}
where instead of the exact data $\bar{r}$, a noisy version $r^{\delta}$ is known. As in \cite{FIM}, if we linearize the forward map $F$ defined in (\ref{fm}) by its Fr\'{e}chet derivative, then by Lemma \ref{Fderivative}, we obtain
$$
F'\theta= F(\theta)-F(0)+o(\|\theta\|),
$$\
i.e.
$$
F'\theta= (u_{\theta}-u)|_{\Gamma_{A}}+o(\|\theta\|).
$$
Here $u_{\theta}|_{\Gamma_{A}}$ and $u|_{\Gamma_{A}}$ are voltages measured in experiments.  In practice they are usually given in  a noisy form as $u_{\theta}|^{\d}_{\Gamma_{A}}$ and $u|^{\d}_{\Gamma_{A}}$ with $\d$ being the noise level for unavoidable experimental error for the measurements of the voltage. When  $\|\theta\|$ is rather small, one can assume that these measurement errors in  \emph{voltage contrast} $(u_{\theta}-u)|_{\Gamma_{A}}$ have the same order of magnitude  as the truncation error $o(\|\theta\|)$. Thus, we take
$$r^{\delta}:= u_{\theta}^{\d}|_{\Gamma_{A}}-u^{\d}|_{\Gamma_{A}}$$
as the noisy right hand term for (\ref{3.1}) and assume that
\begin{equation}\label{3.2}
 \|\bar{r}-r^{\d}\|_{L^2(\Gamma_{A})}\leq\tilde{\delta}=C\delta,
\end{equation}
where a constant $C$ is unknown.

If the Tikhonov regularization is applied to the ill-posed operator equation
$$F'\theta=r^{\delta},$$ then
 the regularized solution solves
\begin{equation}\label{3.3}
(F')^{\ast} F'\theta+\alpha I=(F')^{\ast}r^{\delta}
\end{equation}
where $\a>0$ is the Tikhonov regularization parameter and $I$ is the
 identity operator on space $H^{1}_{0}(\Gamma_{I})$.
It is well known that the solution to (\ref{3.3}) will be the minimizer of the functional
\begin{equation}\label{3.4}
\mathcal{J}(\theta):=\|F'\theta -r^{\delta} \|^2_{L^2(\Gamma_{A})}+\alpha\|\theta\|^2_{H^{1}_{0}(\Gamma_{I})}.
\end{equation}

Here, we assume the exact solution $\theta$ belongs to
the set of source condition
\begin{equation}\label{3.5}
\mathcal{M}_{h}:=\left\{s\in H^{1}_{0}(\Gamma_{I}): s=h((F')^{\ast}F')w, \|w\|\leq 1\right\}
\end{equation}
where $h$ is an {\it `index function'} defined on $[0,\infty)$, which is operator monotone (see \cite{MP1,MP2}) and
satisfies the condition $h(0) = 0$. Moreover, it has been proven that
\begin{equation}\label{3.6}
\sup_{0< \l \leq b} \left|\frac{\a}{\a+\l}\right|h(\l)\leq  h(\a)\q\h{for all}\q \a\in (0, \bar\a]
\end{equation}
and some $\bar\a>0$.

Let us notice that $\mathcal{J}(\theta)$ in \eqref{3.4} is the standard Tikhonov regularization functional, where the penalty term
naturally is imposed in $H^{1}_{0}$-norm. Such a consideration can facilitate the  analysis for the accuracy. Moreover, it is  equivalent to
the Tikhonov regularization functional considered in  \cite{FIM} with a penalty term based on the derivative of the regularized solution.

  The discretization of the regularized problem (\ref{3.3}) is realized by Galerkin method. The Galerkin approximation of Tikhonov-regularization consists in minimizing the above functional $\mathcal{J}(\theta)$
in a finite-dimensional subspace $X_{n}\subset H^{1}_{0}(\Gamma_{I})$.  As usual, in Galerkin scheme, the discretized regularized solution $\theta^{\d}_{\a,n}$ is characterized by the variational
equations
\begin{equation}\label{3.7}
\langle F' \theta^{\d}_{\a,n}-r^{\d}, F' z \rangle +\a  \langle \theta^{\d}_{\a,n}, z \rangle=0, \;\; \forall z\in X_{n},
\end{equation}
or, equivalently,
\begin{equation}\label{3.8}
\theta^{\d}_{\a,n}=\left((F'_{n})^{\ast}F'_{n}+\a I \right)^{-1} (F'_{n})^{\ast}r^{\d},
\end{equation}
where $F'_{n}:=F'P_{n}$, with $P_{n}$ being the projection from $ H^{1}_{0}(\Gamma_{I})$ onto $X_{n}$.

Let $f_{1}$, $f_{2}$, $\ldots$, $f_{n}$ be basis functions of $X_{n}$. If one decomposes $\theta^{\d}_{\a,n}$ into a linear combination of
$f_{1}$, $f_{2}$, $\ldots$, $f_{n}$, i.e. $\theta^{\d}_{\a,n}=\sum_{i=1}^{n}c_{i}f_{i}$, then the coefficient vector $\vec{c}=\{c_{i}\}_{i=1}^{n}$  can be obtained by solving a linear algebraic system,
$$
(M+\a G)\vec{c } = R_{\d},
$$
with the following matrices and vector,
\begin{eqnarray}
  M &:=& \left(\langle F'f_{i}, F'f_{j}\rangle_{L^2(\Gamma_{A})}\right)_{i,j=1}^{n}\nonumber \\
  G &:=& \left(\langle f_{i},f_{j}\rangle_{ H^{1}_{0}(\Gamma_{I})}\right)_{i,j=1}^{n} \label{3.9}\\
  R_{\d} &:=& \left\{\langle F'f_{i},r^{\d}\rangle_{L^2(\Gamma_{A})}\right\}_{i=1}^{n}\nonumber.
\end{eqnarray}
\begin{remark}
The adjoint operator $(F')^{\ast}$ is not involved in the construction of $\theta^{\d}_{\a,n}$. Theoretically, $F'f_{i}$ can be  obtained by solving the boundary value system \eqref{uprime} and deriving the trace on $\Gamma_{A}$, where function $\theta$ is replaced by $f_{i}$, for $i=1,\ldots,n$. Moreover, we do not need each $F'f_{i}$ in an explicit form, but only its products in \eqref{3.9}, which can be computed much more accurately than $F'f_{i}$ itself.
\end{remark}

According to the classical results on Tikhonov regularization for linear ill-posed problem and in view of  (\ref{3.2}) and (\ref{3.6}), it holds that
\begin{eqnarray*}
&&\left\|\theta-\theta_{\alpha,n }^{\d}\right\|_{H^{1}_{0}(\Gamma_{I})}
 \nonumber\\
 &\leq& \left\| \theta-\left((F'_{n})^{\ast}F'_{n}+\a I \right)^{-1} (F'_{n})^{\ast}\bar{r}\right\|_{H^{1}_{0}(\Gamma_{I})}\\
 &&+\left\|\left((F'_{n})^{\ast}F'_{n}+\a I \right)^{-1} (F'_{n})^{\ast}\left(\bar{r}-r^{\d}\right)\right\|_{H^{1}_{0}(\Gamma_{I})}\nonumber\\
 &\leq& \left\| \theta-\left((F'_{n})^{\ast}F'_{n}+\a I \right)^{-1} (F'_{n})^{\ast}\bar{r}\right\|_{H^{1}_{0}(\Gamma_{I})}+\frac{C\d}{2\sqrt{\a}}. \nonumber
\end{eqnarray*}
As in \cite{MP2}, we can estimate the noise free term as follows,
\begin{eqnarray*}
&&\left\| \theta-\left((F'_{n})^{\ast}F'_{n}+\a I \right)^{-1} (F'_{n})^{\ast}\bar{r}\right\|_{H^{1}_{0}(\Gamma_{I})}\nonumber\\
&\leq&\left\| \left(I-\left((F'_{n})^{\ast}F'_{n}+\a I \right)^{-1} (F'_{n})^{\ast}F'_{n}\right)\theta \right\|_{H^{1}_{0}(\Gamma_{I})}\nonumber\\
&&+
\left\| \left((F'_{n})^{\ast}F'_{n}+\a I \right)^{-1} (F'_{n})^{\ast}F'\left(I-P_{n}\right)\theta \right\|_{H^{1}_{0}(\Gamma_{I})}\nonumber\\
&\leq&\left\| \left(I-\left((F'_{n})^{\ast}F'_{n}+\a I \right)^{-1} (F'_{n})^{\ast}F'_{n}\right)h((F'_{n})^{\ast}F_{n}')w \right\|_{H^{1}_{0}(\Gamma_{I})}\nonumber\\
&&+\left\| \left(I-\left((F'_{n})^{\ast}F'_{n}+\a I \right)^{-1} (F'_{n})^{\ast}F'_{n}\right)\left(h((F')^{\ast}F')-h((F'_{n})^{\ast}F_{n}')\right)w \right\|_{H^{1}_{0}(\Gamma_{I})}\nonumber\\
&&+\frac{\left\|F'\left(I-P_{n}\right)\theta \right\|_{L^2(\Gamma_{A})}}{\sqrt{\a}}\nonumber\\
&\leq & C_{1} \left( h(\a)+ h\left(\left\|F'\left(I-P_{n}\right)\right\|^{2}\right)+\frac{\left\|F'\left(I-P_{n}\right) \right\|}{\sqrt{\a}}\right),
\end{eqnarray*}
where the constant $C_{1}$ does not depend on $\a$ and $n$.

In view of the best possible order of accuracy without discretization being $h(\a)+\d/\sqrt{\a}$, the discretization has to be chosen such that
\begin{equation}\label{3.10}
\left\|  F'\left(I-P_{n}\right):  H^{1}_{0}(\Gamma_{I}) \to   L^2(\Gamma_{A})       \right \| \leq \delta.
\end{equation}

Summing up the estimates above, we have the following theorem.
\begin{theorem}
Under assumptions (\ref{3.2}) and (\ref{3.5}), and with  discretization satisfying (\ref{3.10}),
 there holds that
\begin{equation}\label{3.11}
\left\|\theta-\theta_{\alpha,n }^{\d}\right\|_{H^{1}_{0}(\Gamma_{I})} \leq  \bar{K} h(\a)+K\frac{\d}{\sqrt{\a}}
\end{equation}
where the  constants $\bar{K}$ and $K$ do not depend on $\a$ and $\d$.
\end{theorem}

\section{Parameter choice rule based on the balancing principle}
\label{sec4}
In this section, we give a regularization parameter choice rule based on the balancing principle developed in \cite{CaP1,CaPKli,Lu}. The essential idea of this principle is to choose the regularization parameter $\a$ balancing the two parts in error estimate (\ref{3.11}). As a posteriori parameter choice rule, the balancing principle can select regularization parameter in an adaptive way without a priori knowledge of the solution set (\ref{3.5}). That is, the index function $h$ in the bound (\ref{3.11}), which indicates the smoothness of $\theta$ as shown in (\ref{3.5}), does not need to be known. At the same time, it does not require to know the precise noise level, either. In our model problem, constant $C$ in (\ref{3.2}) indicating the precise noise level is unknown, which leads to $K$ in (\ref{3.11}) is also unknown.  A reference noise level $\d$ is sufficient for the performance of the balancing principle. The regularization parameter chosen by the balancing principle leads to a regularized solution with an order-optimal accuracy.

Assume that the projection $P_{n}$ is chosen with $n=n(F',\d)$ such that (\ref{3.10}) is satisfied. Let $\theta_{\a}^{\d}:=\theta_{\a,n(F',\d)}^{\d}$.

We select parameter $\alpha$ from the geometric sequence
$$
\Delta:=\{  \alpha_{n}=\a_{0}q^{n}, \;\; n=0,1,\ldots, N \},
$$
with $q>1$,  sufficiently small $\a_{0}$, and sufficiently large  $N$ such that $\a_{N-1} \leq 1 < \alpha_{N}$.

For any given $K$ , one can choose the parameter from  $\Delta( \alpha_{0}, q, N)$ by the following adaptive strategy,
\begin{equation}\label{4.1}
\alpha(K)=\max\left\{
 \begin{array}{c}\a_{n}\in \Delta:
 \|\theta^{\delta}_{\a_{n}}-\theta^{\d}_{\a_{m}}\|_{H^{1}_{0}(\Gamma_{I})}\leq  K\delta\left( \frac{3}{\sqrt{\a_{n}}}+\frac{1}{\sqrt{\a_{m}}
 }\right),\\
  m=0,1,\ldots,n-1
  \end{array}
   \right\}.
\end{equation}

We further rely on the assumption that  a two-sided estimate
\begin{equation}\label{4.2}
c K \frac{\d}{\sqrt{\a}}\leq \|\theta_{\a}^{0}-\theta_{\alpha}^{\d}\|_{H^{1}_{0}(\Gamma_{I})}\leq  K \frac{\d}{\sqrt{\a}}.
\end{equation}
holds for some $c\in (0,1)$, where $\theta_{\a}^{0}$ is defined by (\ref{3.8}) with $r^{\d}$ being taken as $\bar{r}$. The upper estimate in (\ref{4.2}) is due to (\ref{3.11}). As to the lower estimate, it just suggests that the noise propagation error is not that small.  If the lower estimate is not satisfied, it just means that our estimate to noise level is too pessimistic. However this will not cause a problem, since later we shall show that under assumption (\ref{4.2}) balancing principle can provide an order-optimal accuracy.

Now, consider the following hypothesis set of possible values of the constant $K$
$$
\mathcal{K}=\left\{k_{j}=k_{0}p^{j}, \;\; j=0,1, \ldots, M\right\}, \;\;\;\;p>1,
$$
and assume that  there are two adjacent terms
$k_{l}, k_{l+1}\in\mathcal{K} $  such that
\begin{equation}\label{4.3}
k_{l}\leq cK \leq K \leq k_{l+1}.
\end{equation}

In fact, each element in $\mathcal{K}$ can be viewed as a candidate for the estimator to constant $K$. Our aim is to detect $k_{l+1}$ (or say $k_{l}$) among
the elements in   $\mathcal{K}$, and to use $k_{l+1}$ in adaptive strategy (\ref{4.1}) to obtain a parameter $\a$.

In view of (\ref{4.2}) and (\ref{4.3}), if the hypothesis $k_{j}\in \mathcal{K}$ for $K$ is chosen too small, i.e., $k_{j}\leq k_{l}$  then, as it is shown in \cite{CaP1,Lu},  the
corresponding regularization parameter $\a(k_{j})$ will be smaller than a threshold depending on $\a_{0}$ and $p$. Thus, if
\begin{equation}\label{4.4}
\a(k_{i}):=\min\left\{ \a(k_{j}) \geq 9\a_{0}\left(\frac{p^2+1}{p-1}\right)^2\right\},
\end{equation}
then  there holds that, either $i=l$, or $i=l+1$.

In order to guarantee the regularized solution is stable enough, we choose final regularization parameter as
$$
\a_{+}=\a(k_{i+1}).
$$
With such a choice  $\a_{+}$, we have the following theorem.
\begin{theorem}\label{tbal1}
Under the assumptions above,  there holds
$$
\|\theta-\theta_{\alpha_{+}}^{\d}\|_{H^{1}_{0}(\Gamma_{I})}\leq 6 p^2 \sqrt{q}\bar{K} h (\tilde{h}^{-1}(K\delta)),
$$
where $\tilde{h}(t)=\bar{K}h(t)\sqrt{t}$, $\bar{K}$ is the constant from estimation (\ref{3.11}), and $\tilde{h}^{-1}$ is the inverse function of $\tilde{h}$.
\end{theorem}
Note that from \cite{MP1}  it follows that
the  error bound
indicated in Theorem \ref{tbal1}
is order-optimal, i.e., it is only worse by a constant factor  $3 p^2 \sqrt{q}$  than \emph{a priori}  optimal bound $2\bar{K}h (\tilde{h}^{-1}(K\delta))$.
If  index function $h$ in the source condition (\ref{3.5}) is given as $h(\l)=c\l^{\nu}$, $0<\nu\leq 1$,
then $h (\tilde{h}^{-1}(K\delta))=O(\d^{\frac{2\nu}{2\nu+1}})$, which coincides with the classical rate for Tikhonov regularization.

\begin{remark}
The proof of Theorem \ref{tbal1} can be referred to \cite{CaP1} or \cite{CaPKli}.
 For the general discussions on the application of the balancing principle with two flexible parameters, one can refer to \cite{Lu}.
\end{remark}

\section{Numerical tests}

In this section, we present some numerical examples to illustrate the theoretical results obtained above.

In Examples \ref{Ex1}-\ref{Ex3}, we consider the corrosion problem in (\ref{P}) with
\begin{eqnarray*}
&&\Omega=(0, \pi) \times (0,1),\\
&&\Gamma_{A}=(0, \pi)\times \{0\},\\
&&\Gamma_{I}=(0, \pi)\times \{1\},\\
&&\Gamma_{D}=\{0\}\times(0,1)\cup  \{\pi\}\times (0,1).
\end{eqnarray*}
On such a rectangle domain, if the flux $\Phi=\sin (x)$ is given on $\Gamma_{A}$, then the solution to (\ref{P}) is in the form of
$$
 u(x,y)=\left(-\sinh (y)+\frac{\g\sinh(1)+\cosh(1)}{\sinh(1)+\g\cosh(1)}\cosh y\right)\sin (x).
$$
with $\g>0$. We test the same flux $\Phi=\sin (x)$ in Examples \ref{Ex1}-\ref{Ex3}.
\begin{example}\label{Ex1}
\end{example}
In this example the vector field $\theta$ is given as $\theta=(0,\theta_2(x))$ with
$$
\theta_{2}(x)=\int_{0}^{x}\theta_{2}'(t)dt, \; \mbox{ where } \theta_{2}(0)=0,
$$
and
$$
\theta_{2}'(x)=\left\{
\begin{array}{ll}
\dfrac{-\cot (x)+\g \ds\sqrt{\cot^{2} (x)-\g^2+1}}{\cot^{2} (x)-\gamma},& 0\leq x < \dfrac{\pi}{2},\\
\dfrac{-\cot (x)-\g \ds\sqrt{\cot^{2} (x)-\g^2+1}}{\cot^{2} (x)-\gamma},& \dfrac{\pi}{2} \leq x \leq \pi,
\end{array}
\right.
$$
with the constant $\gamma$ such that $0<\gamma<1$.
Such a choice of $\theta$ corresponds to
$$
u_{\theta}(x,y)=\exp(-y)\sin (x)
$$
solving \eqref{Ptheta}.
\begin{figure}[!htb]
\footnotesize
\def\figwidth{0.4\textwidth}
\advance\tabcolsep by -1.5ex
\centerline{
\begin{tabular}{cc}
\epsfig{file=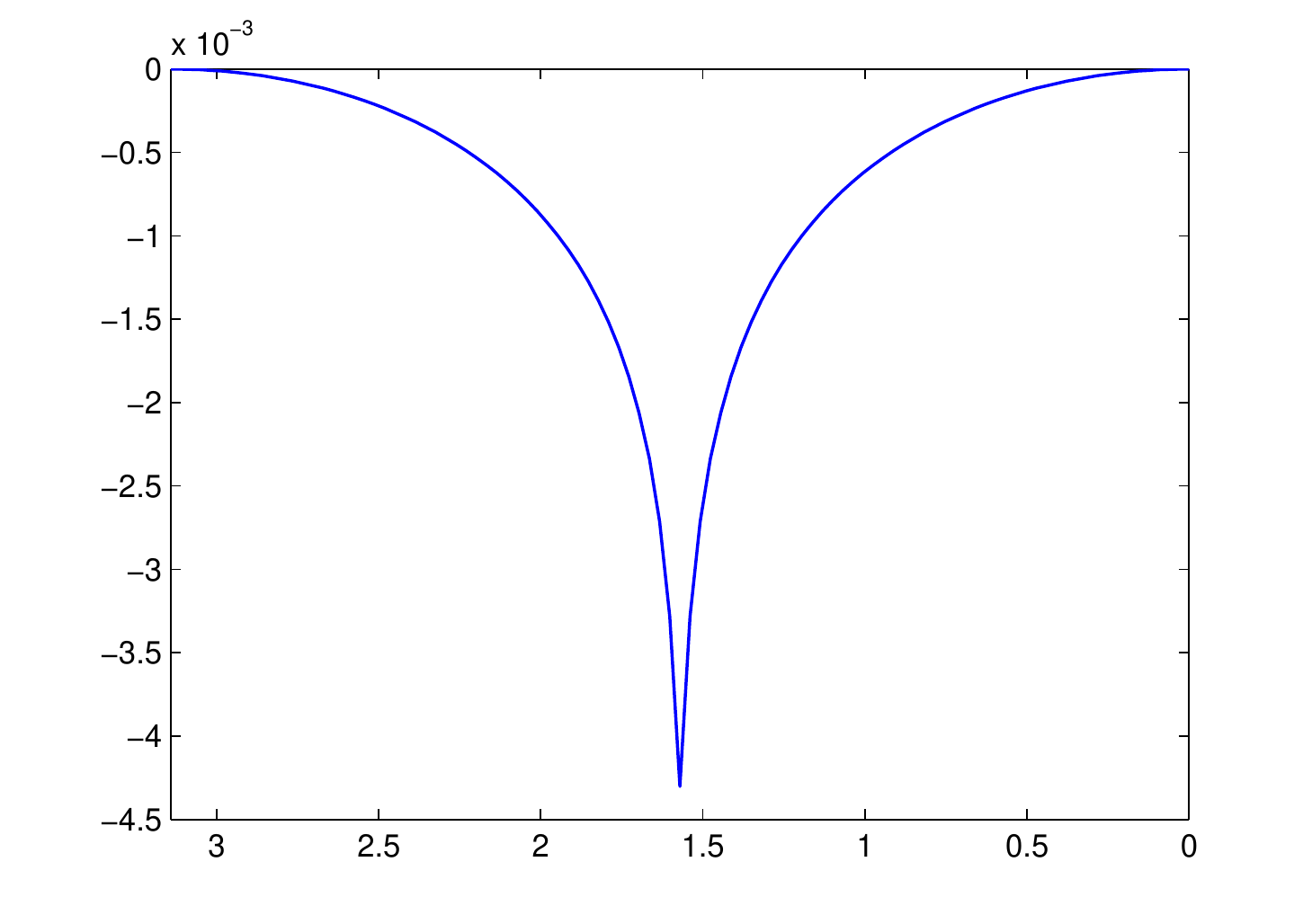,width=\figwidth,height=3cm} &
\epsfig{file=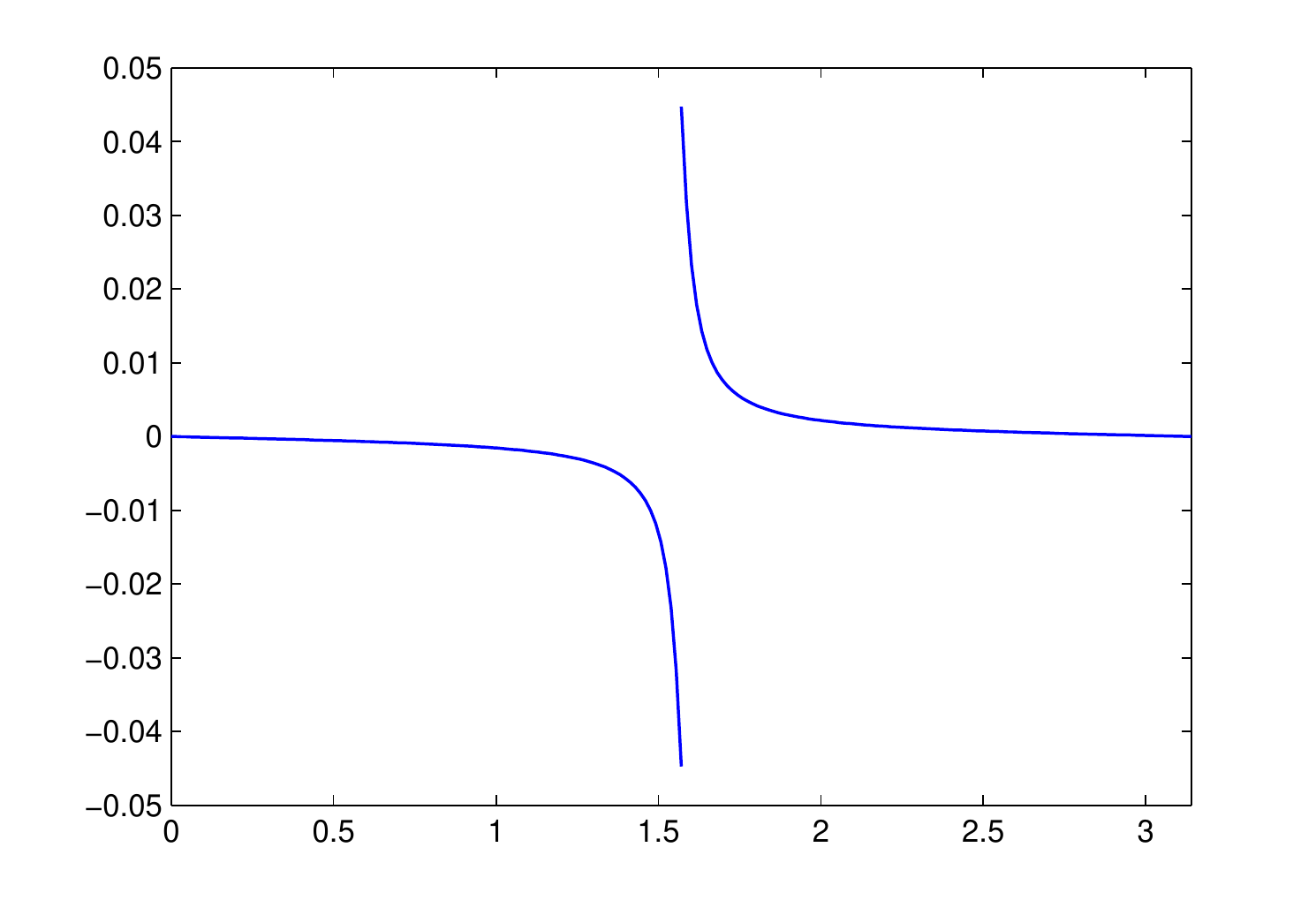,width=\figwidth,height=3cm} \\
(a) $\theta_{2}(x)$ & (b) $\theta_2'(x)$\\
\epsfig{file=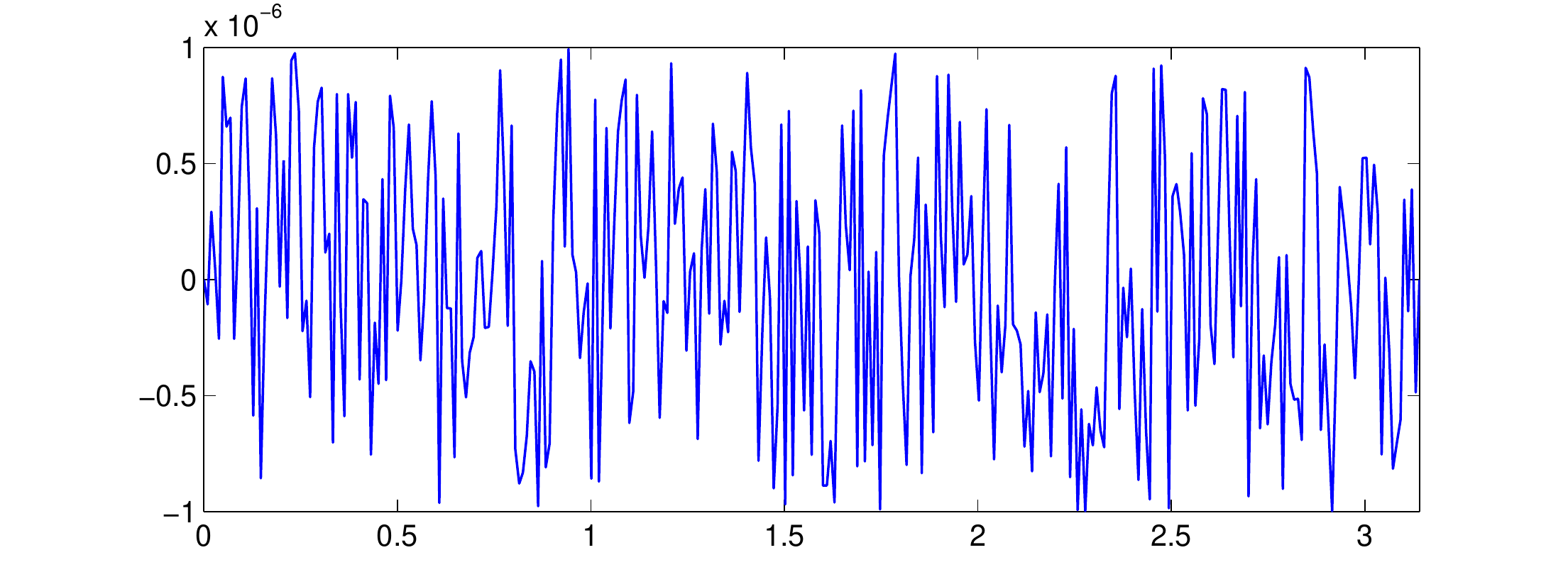,width=\figwidth,height=3cm} &
\epsfig{file=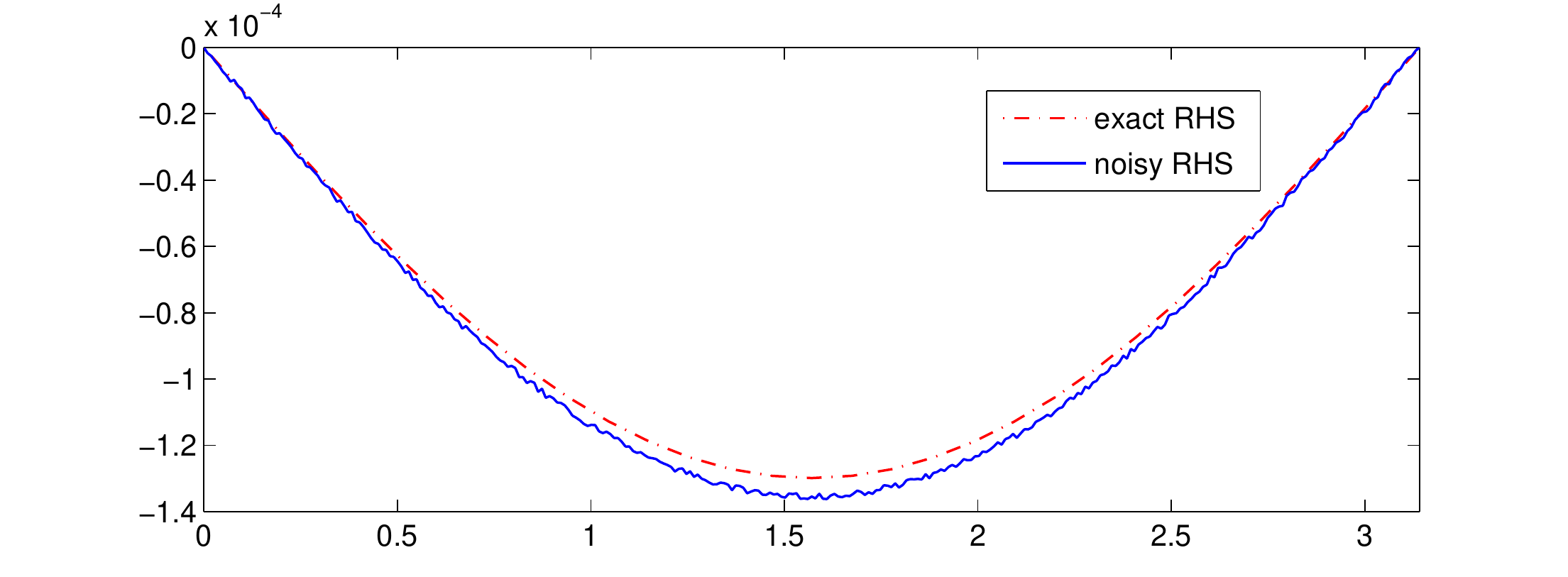,width=\figwidth,height=3cm} \\
\small{(c) noise added to $(u_{\theta}-u)|_{\Gamma_{A}}$  }&
 \small{(d) comparison of $\bar{r}=u^{\prime}|_{\Gamma_{A}}$ and $r^{\d}$}
\end{tabular}}
\caption{\label{fig1}\small{Functions in Example  \ref{Ex1} with $\gamma=0.999$. }}
\end{figure}

As we mentioned in the Introduction, the impedance coefficient $\gamma$ depending on the exterior environment should be a fixed constant in the model problem (\ref{P}). However, in this particular example, the scale of $\theta$ depends on $\gamma$. On the other hand, our linearization approach by truncation can only work when $\|\theta\|$ is rather small. Thus, in this example we test different values of  $\gamma$ which are all quite close to $1$. Figures \ref{fig1} (a) and (b) illustrate the behaviors of $\theta_{2}(x)$ and $\theta_2'(x)$ when $\gamma=0.999$. In order to simulate the error arising in experiment measurements, we add random noise to each grid involved in calculation, i.e. we take
$$
r^{\d}=(u_{\theta}-u)|_{\Gamma_{A}} + \xi \d,
$$
where $\xi$ is random variable with range $[-1,1]$ and the reference noise level $\d=10^{-6}$. Figure \ref{fig1} (c) and (d) show the additional noise $\xi \d$ and the comparison of $\bar{r}=u^{\prime}|_{\Gamma_{A}}$ and $r^{\d}$ in the case $\gamma=0.999$.
\begin{table}[!hbp]
\centering
\begin{tabular}{c||ccccc}
  \hline
  \hline
$\gamma$ &   $\delta$ & $K$ &  $\a_{+}$ & $Err_{h1} $& $Err_{l2}$ \\
\hline
$0.9999$ &   $10^{-8}$ & $0.0489$& $2.015\cdot 10^{-8}$& $0.0012$ & $9.134\cdot 10^{-5}$  \\
\hline
$0.999$ &  $10^{-6}$& $0.0290$& $3.406\cdot10^{-8}$ & $0.0080$ & $ 6.127\cdot 10^{-4}$ \\
\hline
$0.99$ &  $10^{-5}$& $0.0636$& $1.265\cdot10^{-7}$& $0.0438$ & $ 0.0038$  \\
\hline
$0.95$ &  $10^{-4}$& $ 0.0489$& $5.756\cdot10^{-8}$& $ 0.1482$ & $0.0308$  \\
  \hline
  \hline
\end{tabular}
\caption{\label{T1} Test results for Example \ref{Ex1}.}
\end{table}

In all of the following tests,  the discretization level in (\ref{3.7}) is taken as $n=20$ and the regularization parameter $\a$ is chosen by the balancing principle
described in Section \ref{sec4}. In the case of  $\gamma=0.999$, the parameters in the implementation of the balancing principle are settled as follows:
\begin{eqnarray*}
&&\Delta= \{  \alpha_{n}=\a_{0}q^{n}, \;\; n=0,1,\ldots, N \}, \;\;\;\;\a_{0}=1\cdot10^{-11},\;\; q=1.3,\;\; N=69;\\
&&\mathcal{K}=\left\{k_{j}=k_{0}p^{j}, \;\; j=0,1, \ldots, M\right\}, \;\;\;\; k_{0}=0.006,\;\; p=1.3, \;\; M=19.
\end{eqnarray*}
Sequence $\{\a(k_{j})\}_{j=0}^{19}$ produced by \eqref{4.1} with $K$ replaced by $k_{j}$ results in
$$
\begin{array}{ccccc}
 4.827\cdot10^{-11}, &8.157\cdot10^{-11}, &1.379\cdot10^{-10}, &3.937\cdot10^{-10}, &1.900\cdot10^{-9},\\
 9.173\cdot10^{-9},  &3.406\cdot10^{-8},  &7.482\cdot10^{-8},  &9.728\cdot10^{-8}, &1.265\cdot10^{-7},\\
 1.644\cdot10^{-7},  &2.778\cdot10^{-7},  &3.612\cdot10^{-7},  &6.104\cdot10^{-7}, &1.341\cdot10^{-6},\\
 2.946\cdot10^{-6},  &8.415\cdot10^{-6},  &1.094\cdot10^{-5},  &1.422\cdot10^{-5}, &1.849\cdot10^{-5}.
\end{array}
$$
For the parameters  designed as above, the value of the  threshold is calculated as  $9\a_{0}\left(\frac{p^2+1}{p-1}\right)^2= 7.236\cdot10^{-9}$.
Then ${\a}_{+}=\a(k_{7})=3.406\cdot10^{-8}$. At the same time, we obtain an estimate to
 $K$ as $k_{7}=0.0290$, which suggests the true noise level $K\d$ is about $2.90\cdot 10^{-8}$.
\begin{figure}[!htb]
\footnotesize
\def\figwidth{0.45\textwidth}
\advance\tabcolsep by -1.5ex
\centerline{
\begin{tabular}{cc}
\epsfig{file=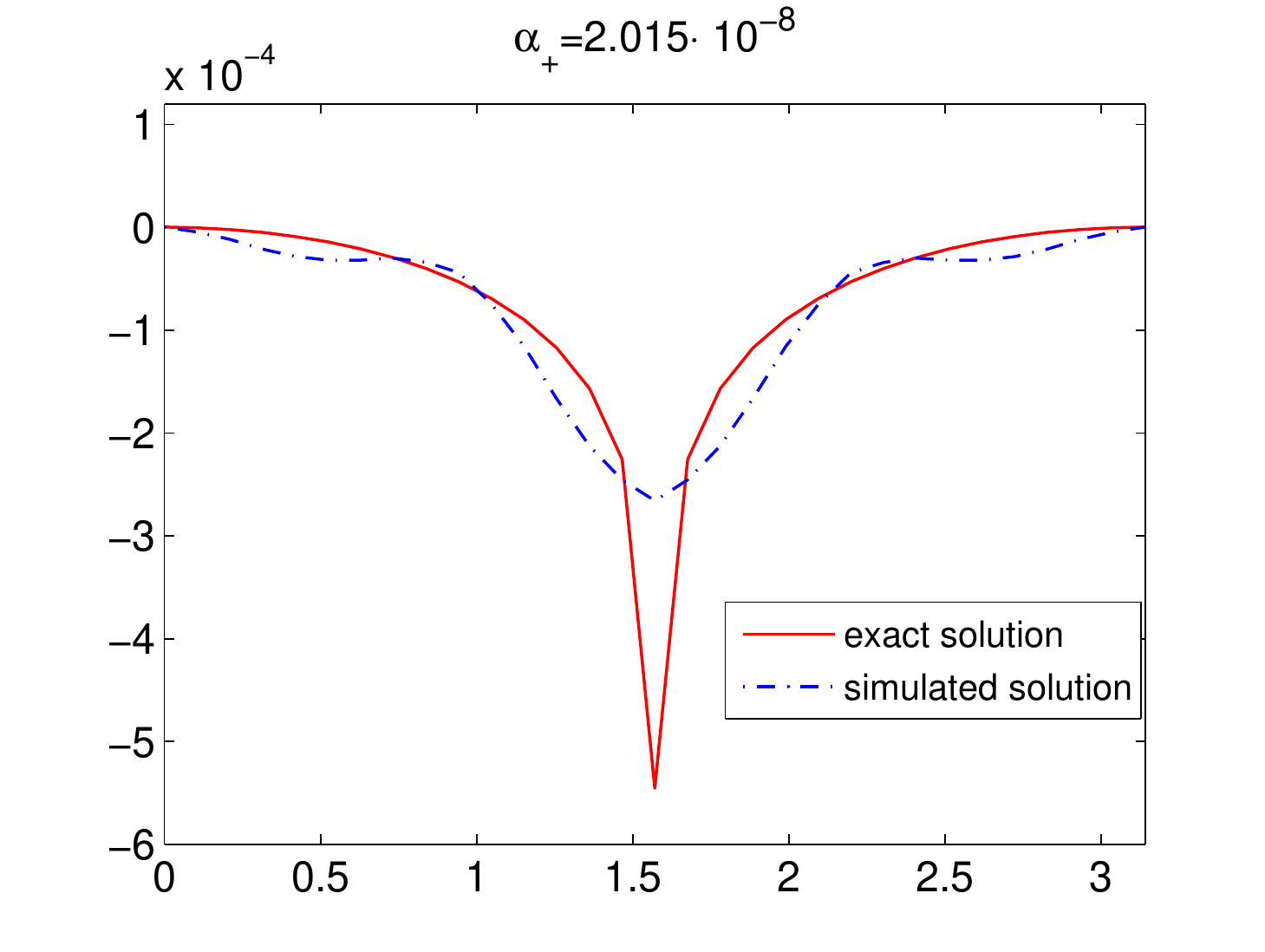,width=\figwidth,height=3.4cm} &
\epsfig{file=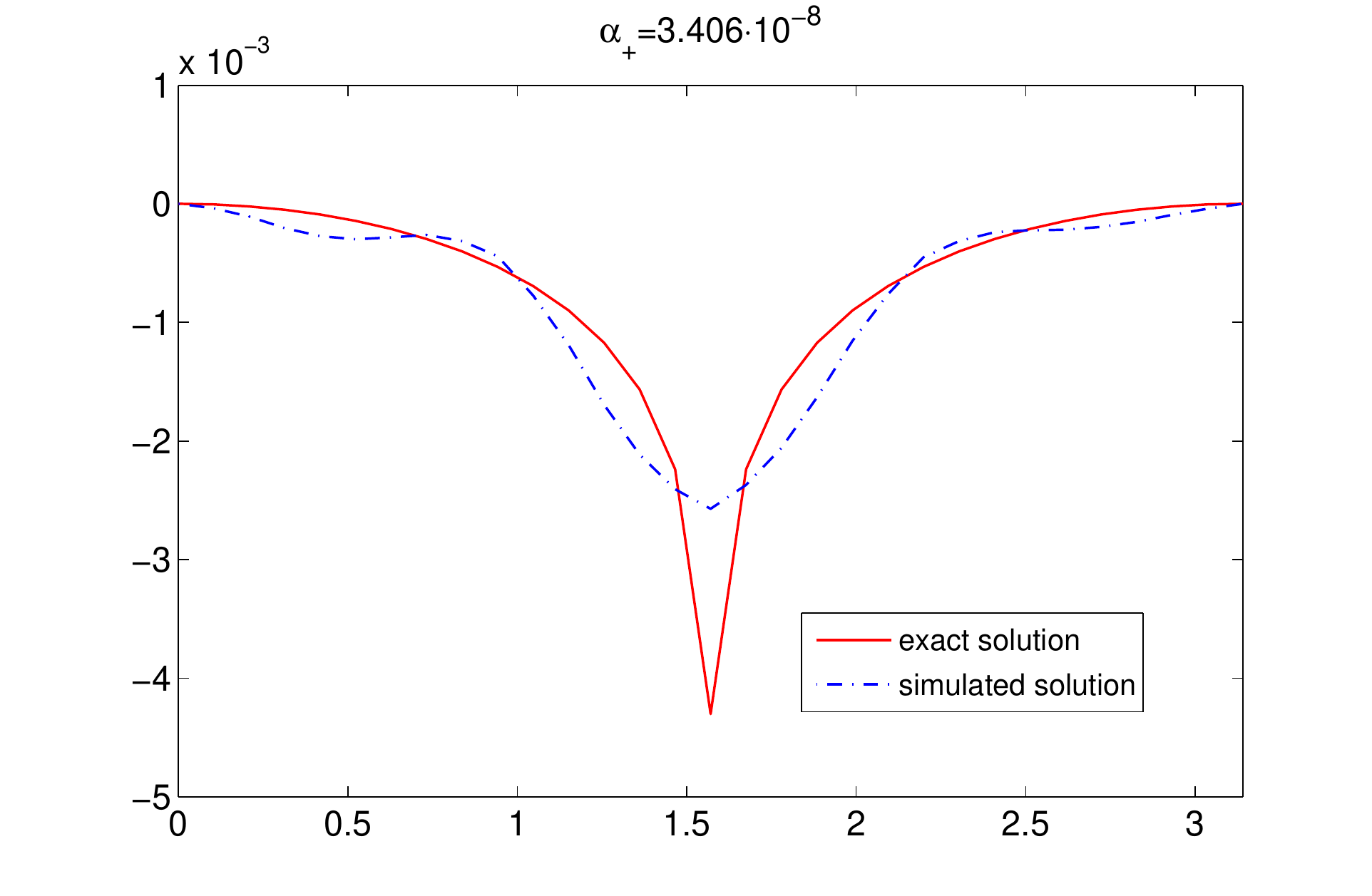,width=\figwidth,height=3.4cm} \\
(a) $\gamma=0.9999$ & (b) $\gamma=0.999$\\
\epsfig{file=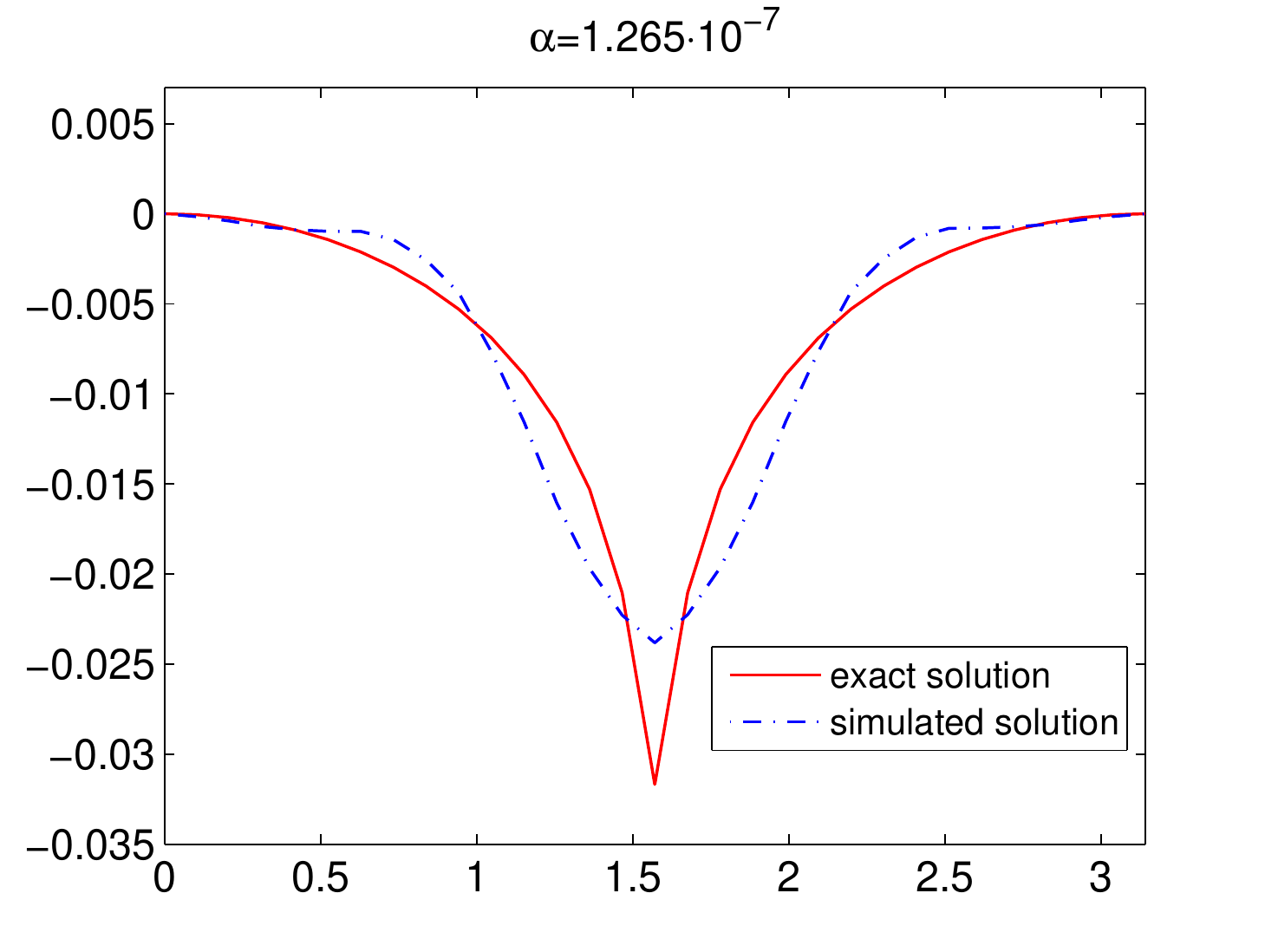,width=\figwidth,height=3.4cm} &
\epsfig{file=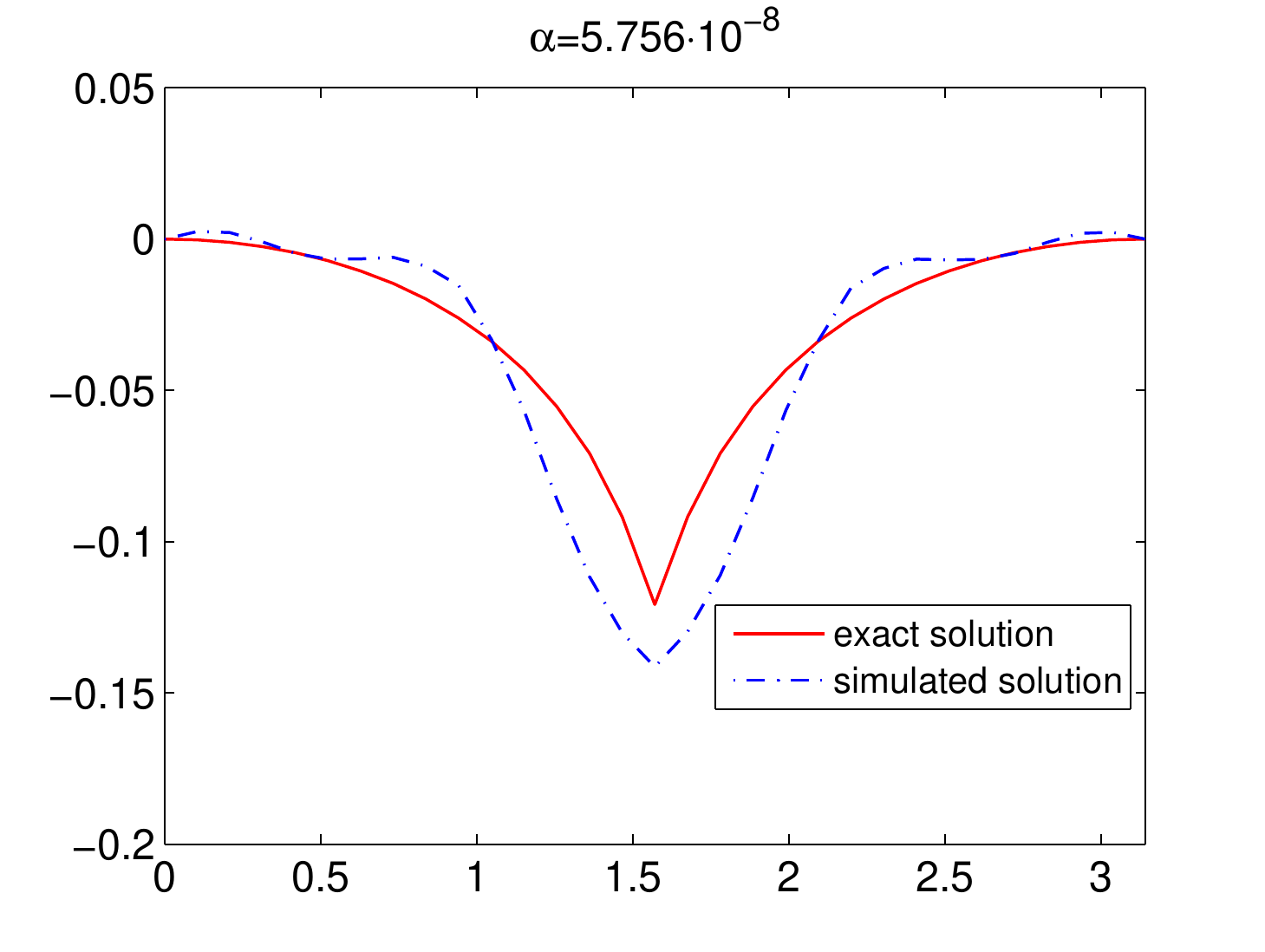,width=\figwidth,height=3.4cm} \\
\small{(c) $\gamma=0.99$  }&
 \small{(d) $\gamma=0.95$}
\end{tabular}}
\caption{\label{fig2}\small{ The simulated solutions $\theta_{\a_{+}}^{\d}$ in Example \ref{Ex1}. }}
\end{figure}

Table \ref{T1} summarizes the results for the other values of $\gamma$ in Example \ref{Ex1}.  Here we take different  reference noise level $\d$ according to $\gamma$, because, as we mentioned, in this particular example, $\gamma$ determines the scale of $\theta$ and furthermore the truncation error. In Table \ref{T1}, $Err_{h1}:=\|\theta-\theta_{\a_{+}}^{\d}\|_{H^{1}(\Gamma_{I})}$ and $Err_{l2}:=\|\theta-\theta_{\a_{+}}^{\d}\|_{L^{2}(\Gamma_{I})}$ denote the errors in the corresponding  norms.
The reconstructed functions $\theta_{\a_{+}}^{\d}$ are displayed in Figure \ref{fig2}.

\begin{example}\label{Ex2}
\end{example}

\begin{figure}[!htb]
\footnotesize
\def\figwidth{0.45\textwidth}
\advance\tabcolsep by -1.5ex
\centerline{
\begin{tabular}{cc}
\epsfig{file=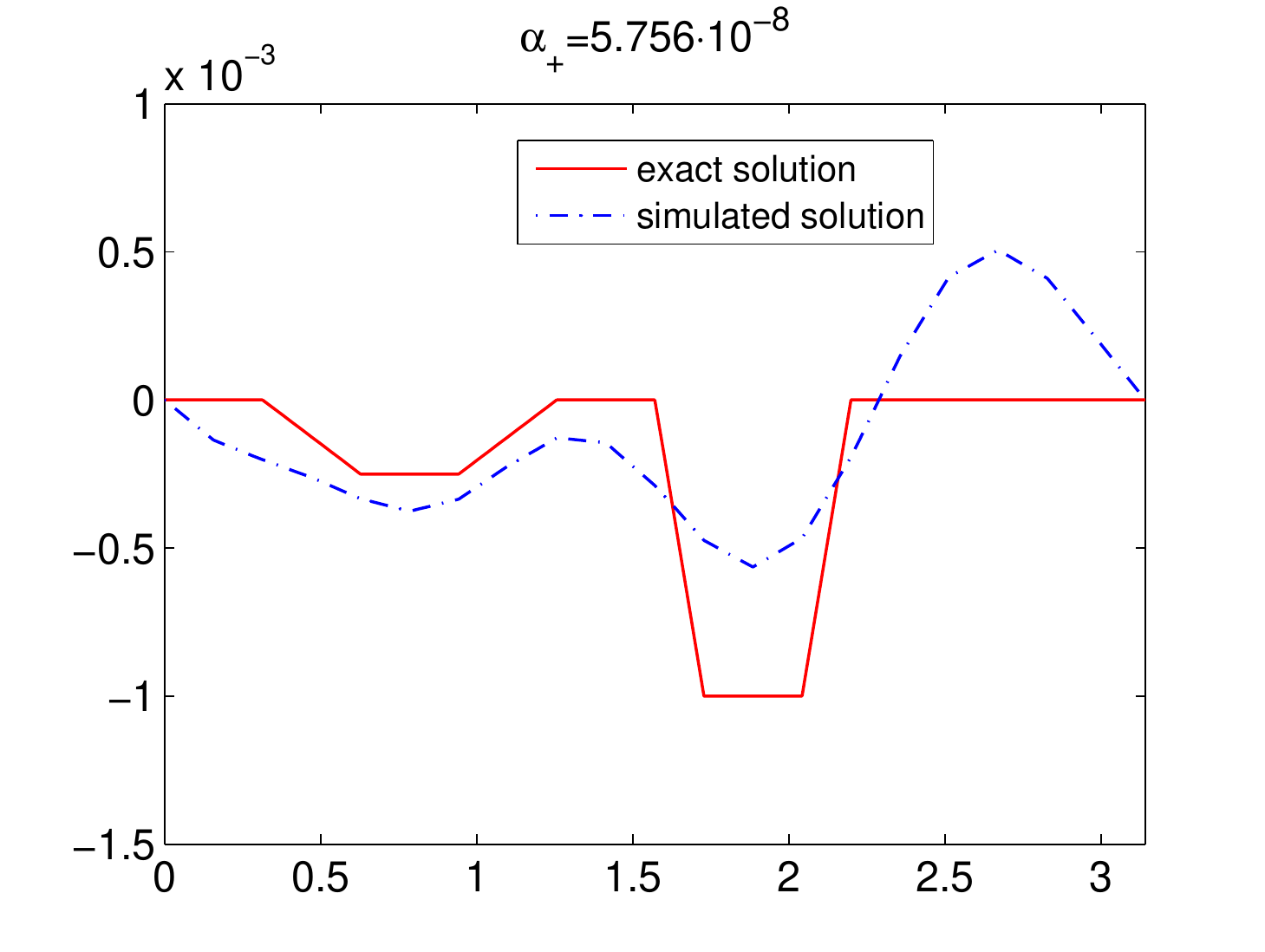,width=\figwidth,height=3.4cm} &
\epsfig{file=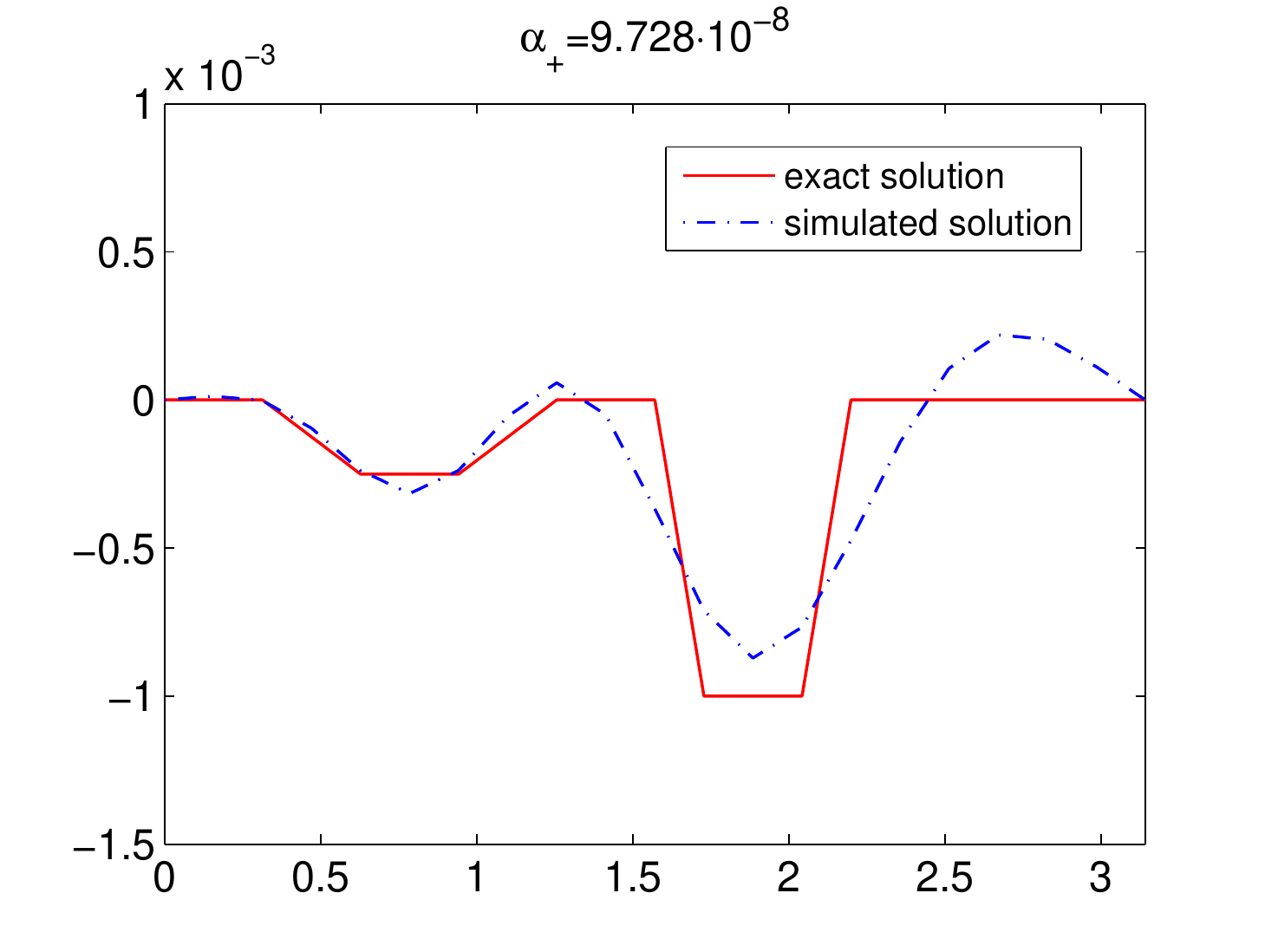,width=\figwidth,height=3.4cm} \\
(a) $\gamma=1$ & (b) $\gamma=10$\\
\end{tabular}}
\caption{\label{fig3}\small{ The simulated solutions $\theta_{\a_{+}}^{\d}$ in Example \ref{Ex2}. }}
\end{figure}
\begin{table}[!htp]
\centering
\begin{tabular}{c||ccccc}
\hline
  \hline
$\gamma$ & $\delta$ & $K$ &  $\a_{+}$ & $Err_{h1} $& $Err_{l2}$ \\
\hline
$1$ & $10^{-7}$ & $0.3937$& $5.756\cdot 10^{-8}$& $ 0.0031$ & $3.565\cdot 10^{-4}$  \\
\hline
$10$  & $10^{-7}$& $0.5119$& $9.728\cdot10^{-8}$ & $0.0030$ & $  2.622\cdot 10^{-4}$ \\
  \hline
\end{tabular}
\caption{\label{T2} Test results for Example \ref{Ex2}.}
\end{table}
In this example, the vector field $\theta=(0,\theta_2(x))$ to be identified is similar to what is considered in \cite{FIM},
where $\theta_{2}(x)$ is a piecewise linear function, as shown in Figure \ref{fig3}. In contrast to Example \ref{Ex1}, we do not assume that $\gamma\in(0,1)$, and test two cases: $\g=1$, $\g=10$. The solution of (\ref{Ptheta}) and its trace $u_{\theta}|_{\Gamma_{A}}$ are generated numerically. We simulate point-wise random noise in each discretization note on $\Gamma_{A}$ with reference level $\d=10^{-7}$, and  $\a_{+}$ is chosen according to the balancing principle under such a value of $\d$. The approximations  $\theta_{\a_{+}}^{\d}$ are displayed in Figure \ref{fig3} and the test results are summarized in Table \ref{T2}.

\begin{figure}[!htb]
\footnotesize
\def\figwidth{0.45\textwidth}
\advance\tabcolsep by -1.5ex
\centerline{
\epsfig{file=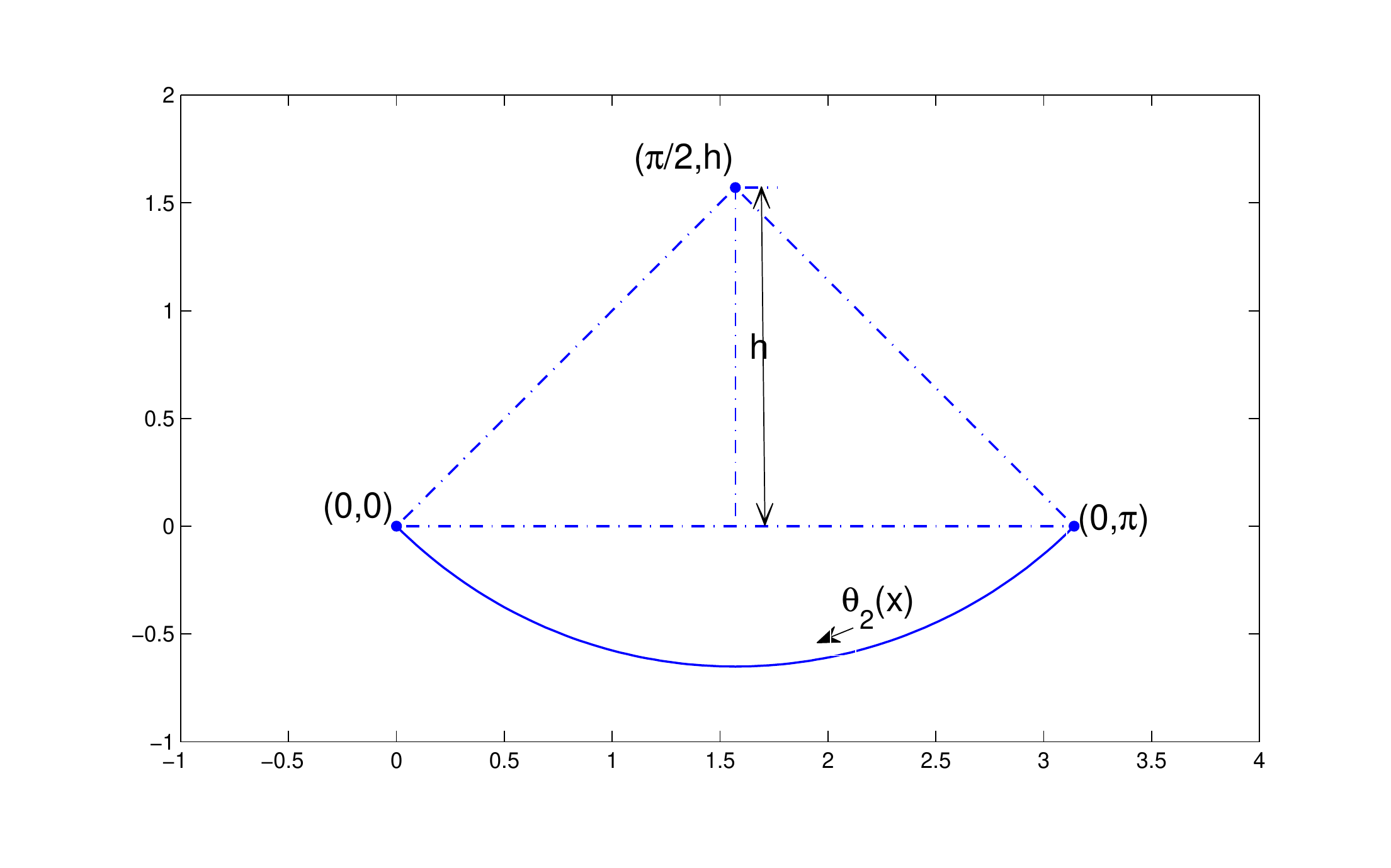,width=\figwidth,height=3.4cm}
}
\caption{\label{fig4}\small{ An illustration for $\theta_{2}(x)$ in Example \ref{Ex3}.}}
\end{figure}

\begin{example}\label{Ex3}
\end{example}
In this example,  we take the vector field $\theta=(0,\theta_2(x))$ with
$$
\theta_{2}(x)=h-\sqrt{\left(\frac{\pi}{2}\right)^2+h^2-\left(x-\frac{\pi}{2}\right)^2},   \;\;\;\; 0< x<\pi,
$$
as  shown in Figure \ref{fig4}.
\begin{figure}[!htb]
\footnotesize
\def\figwidth{0.50\textwidth}
\advance\tabcolsep by -1.5ex
\centerline{
\begin{tabular}{cc}
\epsfig{file=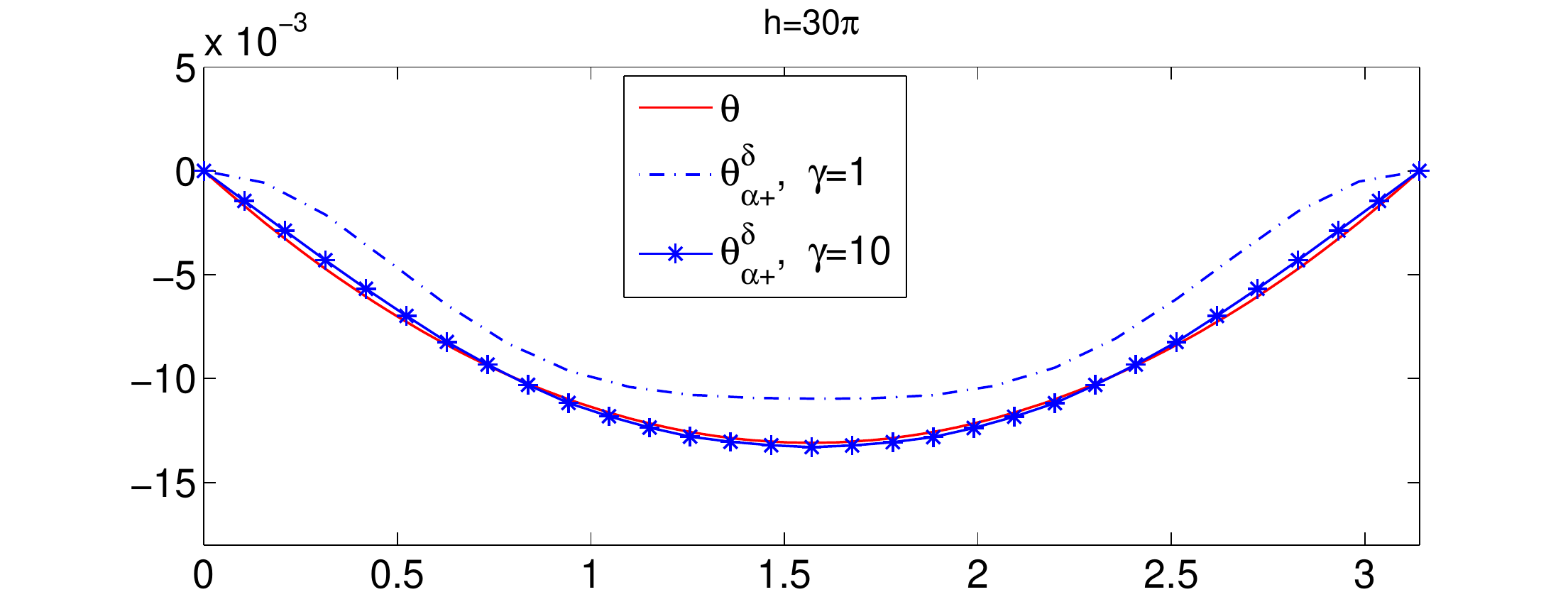,width=\figwidth,height=3.4cm} &
\epsfig{file=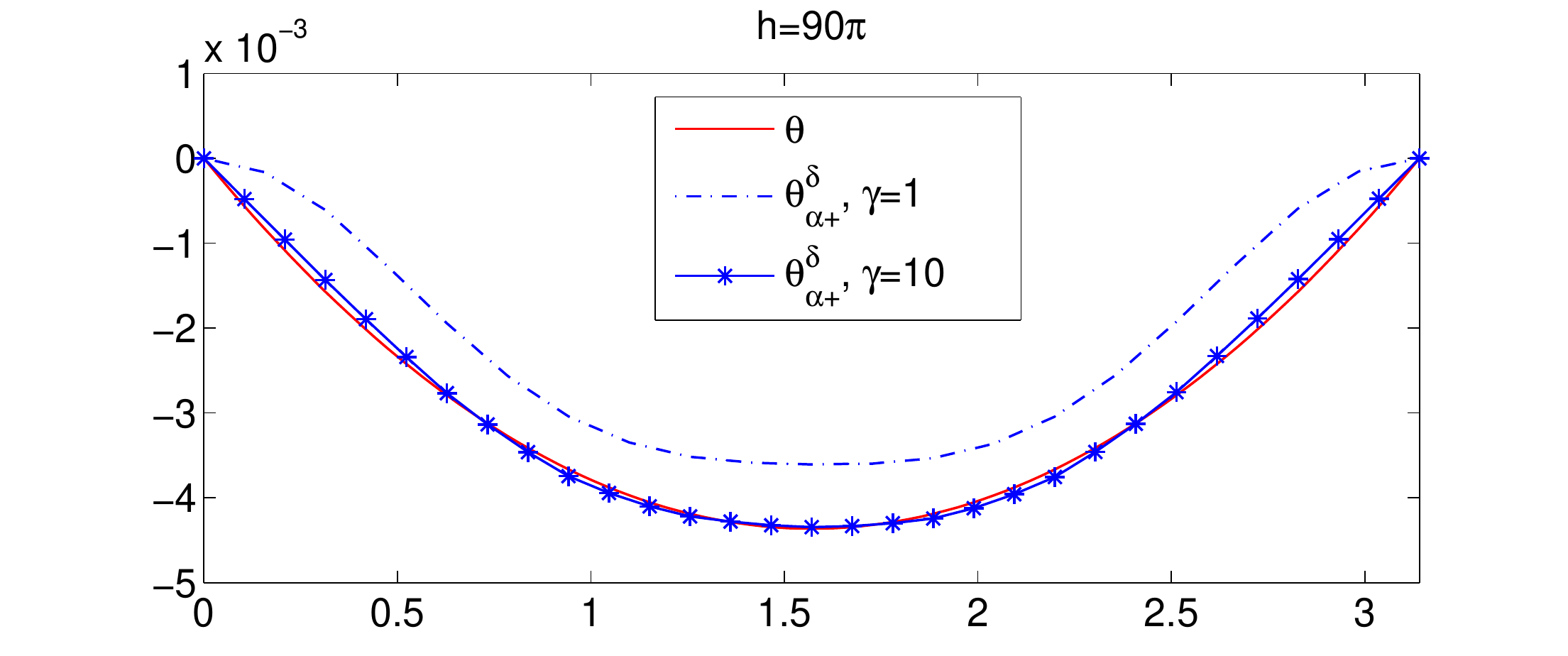,width=\figwidth,height=3.4cm} \\
(a) $h=30\pi$ & (b) $h=90\pi$\\
\end{tabular}}
\caption{\label{fig5}\small{ The simulated solutions $\theta_{\a_{+}}^{\d}$ in Example \ref{Ex3}. }}
\end{figure}
\begin{table}[!htbp]
\centering
\begin{tabular}{ccccccc}
\hline
 \hline
h& $\gamma$ & $\delta$ & $K$ &  $\a_{+}$ & $Err_{h1} $& $Err_{l2}$ \\
\hline
 \multirow{2}{*}{$30\pi$}&  $1$ & $10^{-5}$ & $0.0371$& $1.341\cdot 10^{-6}$& $ 0.0079$ & $0.0034$  \\
\cline{2-7}
 &$10$ & $10^{-5}$& $0.0371$& $1.849\cdot10^{-5}$ & $0.0016$ & $ 4.273\cdot 10^{-4}$ \\
  \hline
  \hline

 \multirow{2}{*}{$90\pi$}  &$1$ & $10^{-7}$ & $0.1379$& $3.611\cdot 10^{-6}$& $  0.0027$ & $0.0013$  \\
\cline{2-7}
& $10$ & $10^{-7}$& $ 0.1060$& $1.743\cdot10^{-6}$ & $6.300\cdot 10^{-4}$ & $ 1.340\cdot 10^{-4}$ \\
  \hline
\end{tabular}
\caption{\label{T3} Test results for Example \ref{Ex3}.}
\end{table}
Here one can change the value of the constant $h>0$ to adjust the scale of $\theta$. The solution $u_{\theta}$ and its trace on $\Gamma_{A}$ in this example are also obtained numerically. In order to guarantee that the scale of $\theta$ is small enough such that the truncation method can work well,  we test  $h=30\pi$ and $h=90\pi$.
In both cases, larger value of  $\gamma$ may make the problem less ill-posed and result in better reconstruction. The approximations  $\theta_{\a_{+}}^{\d}$ are displayed in Figure \ref{fig5} and the test results are summarized in Table \ref{T3}.

\begin{example}\label{Ex4}
\end{example}
\begin{figure}[!htb]
\footnotesize
\def\figwidth{0.4\textwidth}
\advance\tabcolsep by -1.5ex
\centerline{
\epsfig{file=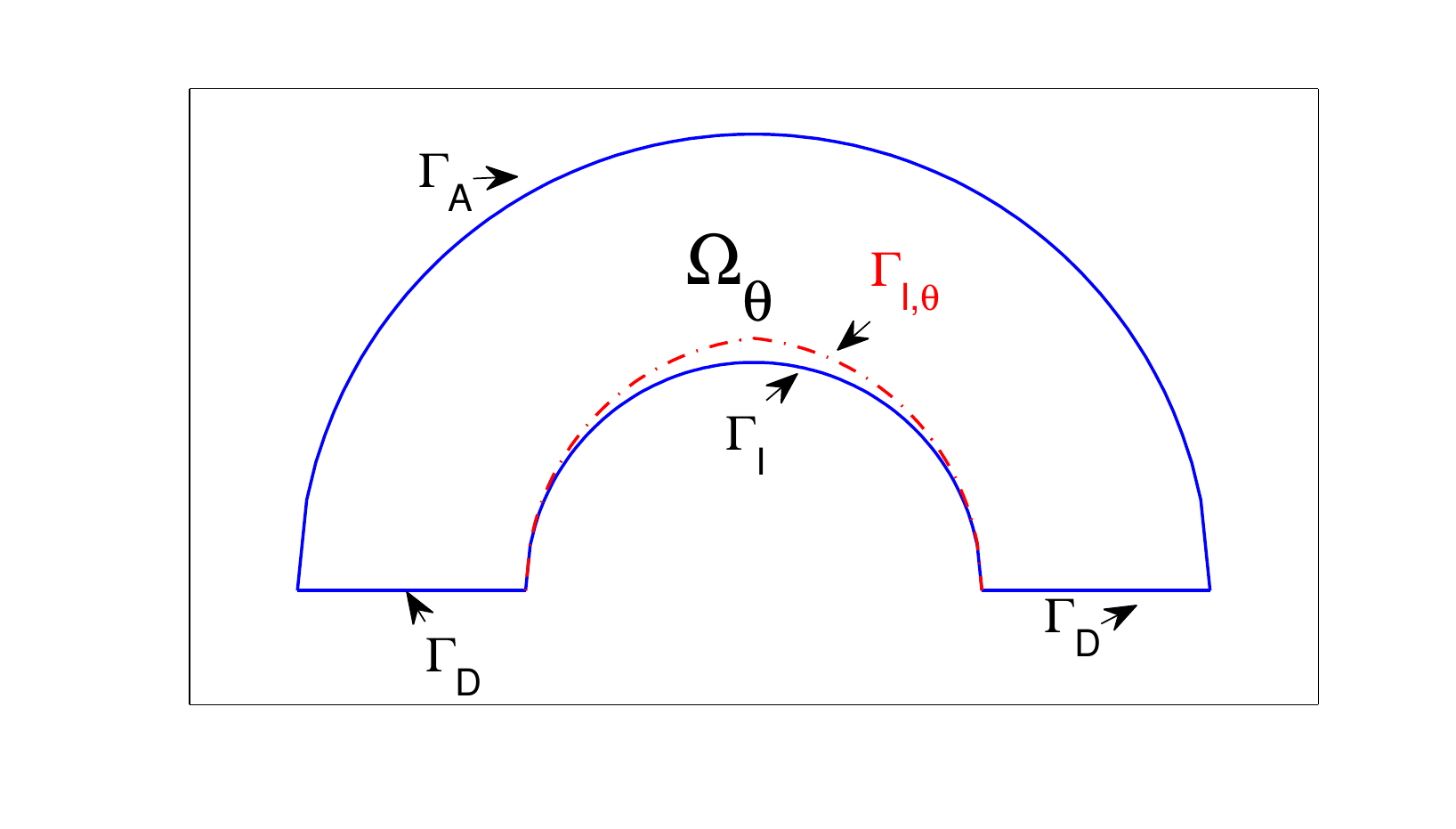,width=\figwidth,height=3cm}
}
\caption{\label{fig6}\small{ Domain $\O_{\theta}$  in Example \ref{Ex4}.}}
\end{figure}
In the last example, we consider a domain $\O$  given as a half annulus bounded by the following curves (see Figure \ref{fig6}).
\begin{eqnarray*}
&&\Gamma_{A}=\left\{(x,y): \;y=\sqrt{4-x^2}, \; -2<x<2\right\},\\
&&\Gamma_{I}=\left\{(x,y): \;y=\sqrt{1-x^2}, \;-1<x<1\right\},\\
&&\Gamma_{D}=\left\{(x,y):\; y=0, \; 1< |x| <2 \right\}.
\end{eqnarray*}
\begin{table}[!htb]
\centering
\begin{tabular}{c||ccccc}
  \hline
  \hline
$\gamma$ &   $\delta$ & $K$ &  $\a_{+}$ & $Err_{h1} $& $Err_{l2}$ \\
\hline
$0.99$ &   $10^{-5}$ & $0.0816$& $1.032\cdot10^{-5}$& $ 0.0057$ & $0.0028$  \\
\hline
$0.9$ &  $10^{-4}$& $0.0816$& $1.032\cdot10^{-5}$ & $0.0708$ & $0.0384$ \\
  \hline
  \hline
\end{tabular}
\caption{\label{T4} Test results for Example \ref{Ex4}.}
\end{table}
\begin{figure}[!htb]
\footnotesize
\def\figwidth{0.45\textwidth}
\advance\tabcolsep by -1.5ex
\centerline{
\begin{tabular}{cc}
\epsfig{file=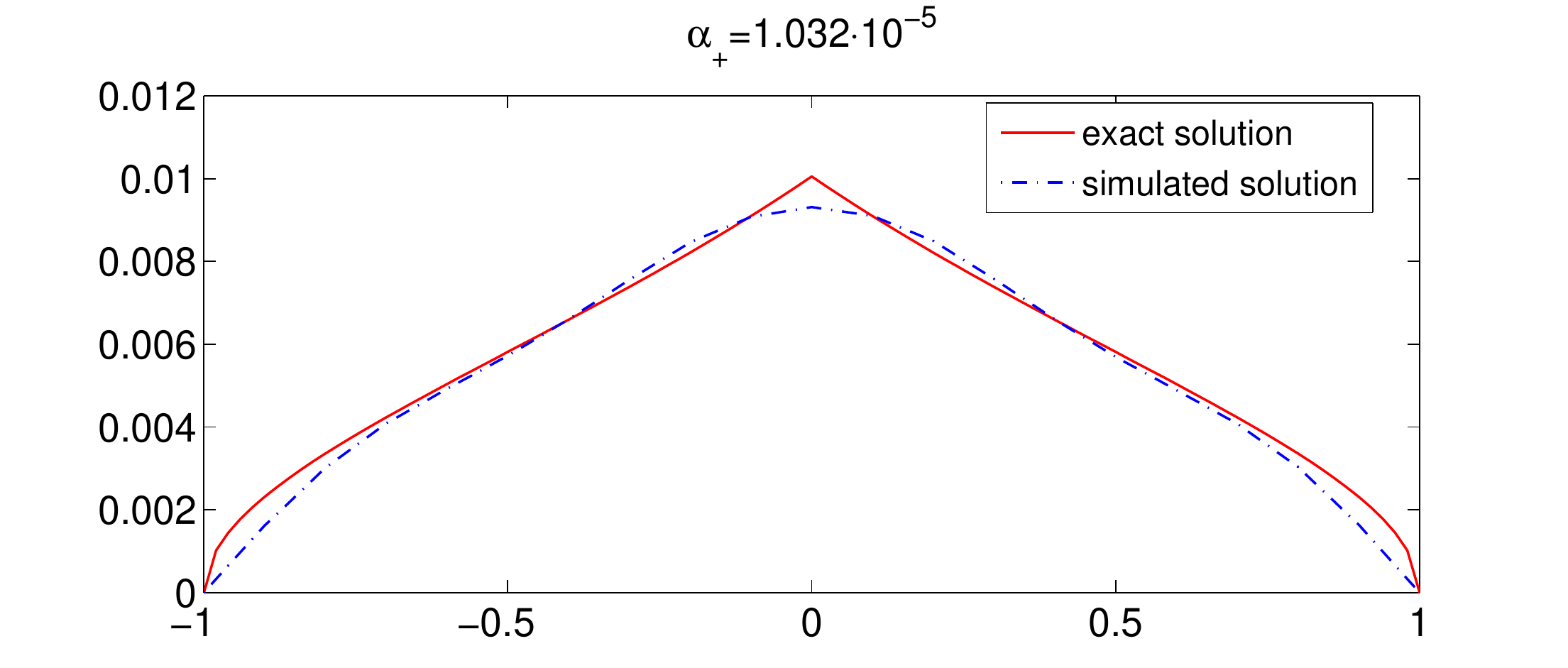,width=\figwidth,height=3.3cm} &
\epsfig{file=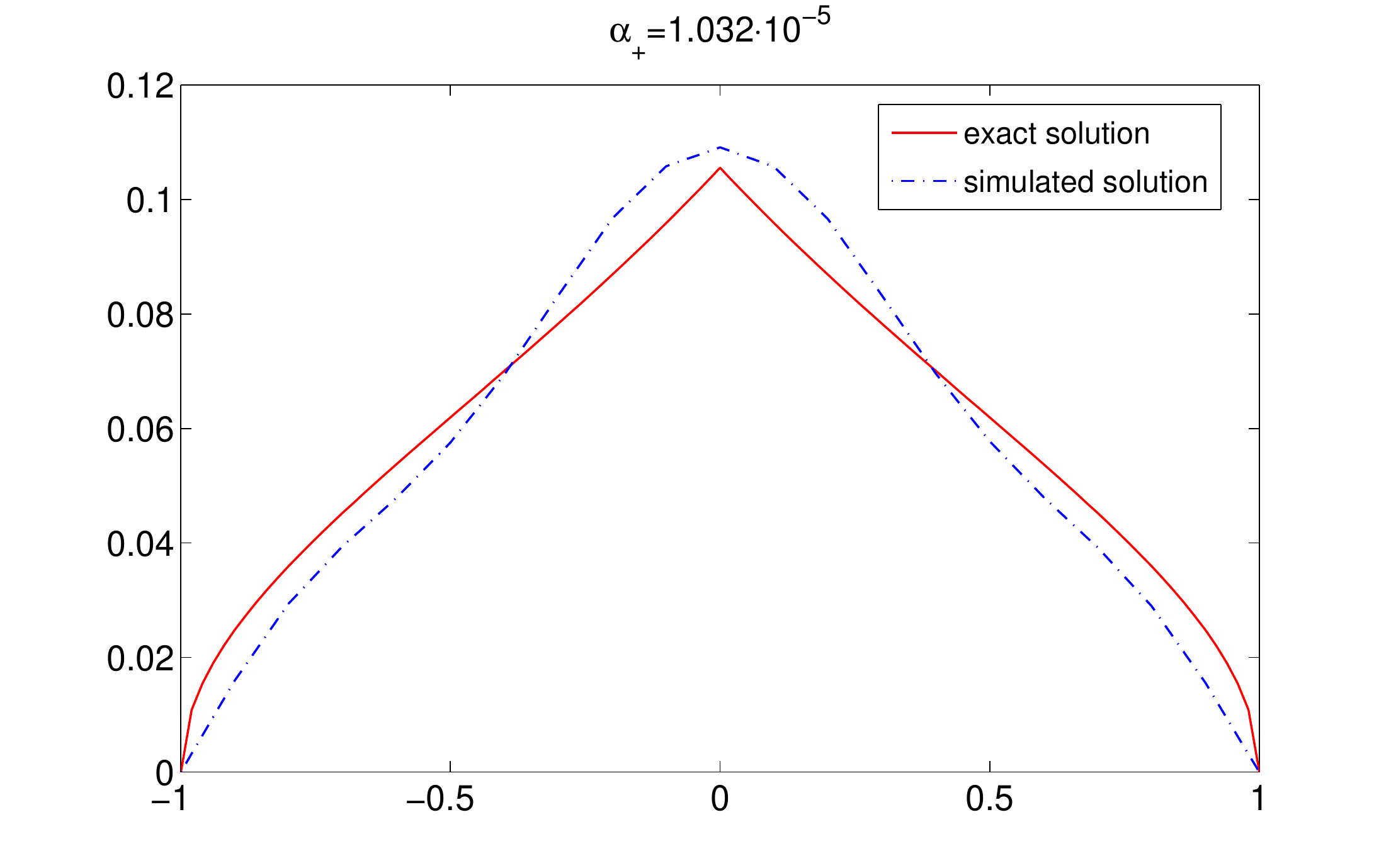,width=\figwidth,height=3.2cm} \\
(a) $\gamma=0.99$ & (b) $\gamma=0.9$
\end{tabular}}
\caption{\label{fig7}\small{ The simulated solutions $\theta_{\a_{+}}^{\d}$ in Example \ref{Ex4}. }}
\end{figure}
For flux $\Phi=y$ on $\Gamma_{A}$, the function $u$ solving (\ref{P}) can be written as
$$
u(x,y)=Ay+\frac{By}{x^2+y^2}, \mbox{ where } A=1+\frac{B}{4} \mbox{ and } B=\frac{1-\gamma}{\frac{5}{4}\gamma+\frac{3}{4}},
$$
with  $0<\g<1$. For vector field $\theta=(0,\theta_2(x))$, we consider $ \theta_{2}(x)=\varphi(x)-\sqrt{1-x^2}$,
where
$$
\varphi(x)=\left\{
\begin{array}{ll}
\dfrac{1}{\gamma}\ds\sqrt{1-(\gamma x+\g-1)^2}, & -1<x \leq 0,\\
\dfrac{1}{\gamma}\ds\sqrt{1-(\gamma x-\g+1)^2}, & 0\leq x \leq 1.
\end{array}
\right.
$$
It can be verified that $u_{\theta}(x,y)=y$
solves\eqref{Ptheta} in $\O_{\theta}$.

 In this example $\gamma$ also determines the scale of $\theta$. Thus, we take the values of $\gamma$ very close to $1$. The approximations  $\theta_{\a_{+}}^{\d}$ are displayed in Figure \ref{fig7} and the test results are summarized in Table \ref{T4}.

We would like to note that in all considered examples the balancing principle (\ref{4.1}), (\ref{4.4}) has been implemented with the same values of the design parameters $\a_{0}$,  $p$ and $q$, because the domain $\O$ and the operator $F'$ are the same for all examples. This suggests that in practice, for a given domain $\O$ the parameters $\a_{0}$,  $p$ and $q$ can be determined in the experiments with a problem (\ref{uprime}) where a solution is known, and then kept for studying all other problems (\ref{Ptheta}) in the given domain $\O$.

Both the theoretical and numerical results suggest that the linearizaton approach considered in this paper can perform well for the identification of the corroded boundary only on condition that the scale of this  boundary function is quite small. This is the limitation of the approach. However, in practice one certainly does not expect too much corrosion taking place to the metallic specimen.

\section*{Acknowledgement}

This research was partially supported by  Natural Science Foundation of China (No.11201497), Foundation for Distinguished Young Talents in Higher Education of Guangdong, China (No. LYM11007), Guangdong Provincial Government
of China through the ``Computational Science Innovative Research Team"
program, Guangdong Province Key Laboratory of Computational Science
at the Sun Yat-sen University and the PRIN 20089PWTPS project funded by the Italian Ministry for Education, University and Research (MIUR). S.V.P. is partially supported by the grant P25424 of the Austrian Science Fund (FWF). E.S. wishes to thank the Department of Mathematics and Geosciences of the University of Trieste (Italy)  for partly supporting her work by a research grant. Part of this work was done while E.S. was visiting the Sun
Yat-sen University. She also wishes to express her gratitude to the School of Mathematics and
Computational Science for the kind hospitality.

\end{document}